\documentclass{amsart}
\usepackage{amssymb,latexsym,amsmath,amsfonts,hhline}
\usepackage{pdfsync}
\usepackage{xcolor} 
\usepackage{comment,colortbl}

\newcommand{\cal}{\mathcal}

\makeatletter
\renewcommand{\subsection}{\@startsection{subsection}{2}{0mm}{-2mm}{-2mm}{\bf\normalsize}}

\makeatother

\newtheorem{formula}{}[section]
\newtheorem{definition}[formula]{Definition}
\newtheorem{corollary}[formula]{Corollary}
\newtheorem{remark}[formula]{Remark}
\newtheorem{lemma}[formula]{Lemma}
\newtheorem{theorem}[formula]{Theorem}

\def\thrm{\begin{theorem}}
\def\thrml#1{\begin{theorem}\label{#1}}
\def\ethrm{\end{theorem}}
\def\rmrk{\begin{remark}}
\def\rmrkl#1{\begin{remark}\label{#1}}
\def\ermrk{\end{remark}}
\def\dfntn{\begin{definition}}
\def\dfntnl#1{\begin{definition}\label{#1}}
\def\edfntn{\end{definition}}
\def\nmrt{\begin{enumerate}}
\def\enmrt{\end{enumerate}}
\def\tm#1{\item[{\rm (#1)}]}
\def\qtn{\begin{equation}}
\def\qtnl#1{\begin{equation}\label{#1}}
\def\eqtn{\end{equation}}
\def\lmm{\begin{lemma}}
\def\lmml#1{\begin{lemma}\label{#1}}
\def\elmm{\end{lemma}}
\def\crllr{\begin{corollary}}
\def\crllrl#1{\begin{corollary}\label{#1}}
\def\ecrllr{\end{corollary}}
\def\css{\begin{cases}}
\def\ecss{\end{cases}}

\def\proof{\noindent{\bf Proof}.\ }

\def\cH{{\cal H}}

\def\cR{{\cal R}}

\def\bN{{\mathsf N}}

\def\fD{{\frak D}}

\def\fK{{\frak K}}

\DeclareMathOperator{\aut}{Aut}

\DeclareMathOperator{\AGL}{AGL}

\DeclareMathOperator{\cyc}{cyc}

\DeclareMathOperator{\iso}{Iso}

\DeclareMathOperator{\poly}{poly}

\DeclareMathOperator{\Sup}{Sup}

\DeclareMathOperator{\sym}{Sym}
\DeclareMathOperator{\Wr}{Wr}

\def\eprf{\hfill$\square$}

\def\qaq{\quad\text{and}\quad}

\def\ov{\overline}

\def\wh{\widehat}

\begin{document}
\title[Testing isomorphism of circulant objects]{Testing isomorphism of circulant objects\\ in polynomial time}
\author{Mikhail Muzychuk}
\address{Ben-Gurion University, Beer-Sheva, Israel,}
\email{muzychuk@bgu.ac.il }
\author{Ilia Ponomarenko}
\address{Steklov Institute of Mathematics at St. Petersburg, Russia}
\email{inp@pdmi.ras.ru}
\thanks{The first author was supported by the Israeli Ministry of Absorption. The second author was supported by the RFBR grant No. 18-01-00752}
\date{}

\begin{abstract}
Let $\fK$ be a class of combinatorial objects invariant with respect to a given regular cyclic group. It is proved that the  isomorphism of any two objects $X,Y\in\fK$ can be tested in polynomial time in sizes of~$X$ and~$Y$.
\end{abstract}

\date{}
\maketitle

\section{Introduction}
There are different ways to define a combinatorial object over a given point set. They include: concrete categories \cite{Babai1977a}, relational structures \cite{P87}, a boolean tower \cite{B} and  hereditarily finite sets \cite{SW}. In practice, every ``reasonable" class of combinatorial objects can be presented using any of them. To formulate our main results, it is convenient to consider combinatorial objects as the objects of   a concrete category or as relational structures.\medskip

In a  {\it concrete category} $\fK$, each object $X\in\fK$ is associated with an underlying set~$\Omega(X)$, and each isomorphism from $X$ to $Y$ is associated with a certain bijection $f:\Omega(X)\to\Omega(Y)$; the set of all these bijections is denoted by $\iso_\fK(X,Y)$. It is also assumed that for any bijection $f$ from the set $\Omega(X)$ to another set, there exists a unique object $Y=X^f$ for which this set is the underlying one and $f\in\iso(X,Y)$. Thus,
$$
X\cong_\fK Y\quad\Leftrightarrow\quad Y=X^f\ \,\text{for some}\ \,f\in\iso(X,Y).
$$
Given a set $K\subseteq\sym(\Omega)$ and two objects $X,Y\in\fK$ with $\Omega(X)=\Omega=\Omega(Y)$, we write $\iso_K(X,Y)$ for the intersection $K\cap \iso_\fK(X,Y)$.\medskip

In what follows, under a Cayley object of $\fK$ over a group $G$, we mean any $X\in\fK$ such that 
$$
\Omega(X)=G\qaq \aut_\fK(X)\ge G_R,
$$
where $\aut_\fK(X)=\iso_\fK(X,X)$ and $G_R$ is  the group induced by the right regular representation of~$G$. In the case of $G$ being cyclic the object will be called {\it cyclic} or {\it circulant}.\medskip

A particular example of a concrete category is formed by relational structures. A {\it relational structure} over a ground set $\Omega$ is a pair $X=(\Omega,\cR)$, where $\cR$ is  a finite set of relations  on~$\Omega$. It is assumed that $\cR$ is linear ordered and the arities of  the relations in $\cR$ may be  different.\footnote{The arity of a relation $R$ on $\Omega$ is defined to the smallest positive integer $\ell$ such that $R\subseteq\Omega^\ell$.} Isomorphisms and automorphisms of relational structures respect the ordering and are defined in a natural way. When $\cR=\{R_1,...,R_n\}$, the number $|\Omega|+\sum_i |R_i|$ is called the {\it size} of~$X$.

The aim of the present note is to provide a complete solution of the following problem.\medskip

{\bf Circulant Objects Isomorphism.} {\it Given a cyclic group $C$ and two Cayley relational structures over~$C$, test whether they are isomorphic and (if so) find an isomorphism between them}.\medskip

The first results about this problem were obtained by Bays and Lambossy \cite{B23,L31}. The group-theoretical approach to the isomorphism problem for Cayley objects was developed by Babai in \cite{Babai1977a}. The first breakthrough towards a solution of the above problem was done by P{\' a}lfy~\cite{P87}. He proved that if the group order $n=|C|$ satisfies the condition $(n,\varphi(n))=1$, then for any two circulant relational structures  $X$ and~$Y$,
$$
\iso(X,Y)\neq \varnothing\quad\Leftrightarrow\quad\iso_{\aut(C)}(X,Y)\neq\varnothing.
$$ 
This result provides a simple polynomial-time algorithm for isomorphism testing of  circulant combinatorial structures of special orders. In order to cover the remaining orders of circulant objects it was proposed in \cite{HJP,H} to replace $\aut(C)$ by a bigger set $S\subset \sym(C)$ with the property 
$$
\iso(X,Y)\neq \varnothing\quad\Leftrightarrow\quad \iso_{S}(X,Y)\neq\varnothing.
$$ 
This idea was further developed in \cite{Muz1999} where such a set was called a {\it solving set}. It was shown in \cite{Muz2004, Muz2011, KKMM} that various classes of circulant combinatorial objects admit solving sets of polynomial size. 
The first main result of our paper shows that there exists a solving set which works for all circulant combinatorial objects.

\thrml{081117a} 
Let $C$ be a cyclic group of order $n$. Then in time $\poly(n)$, one can construct a solvable group $K\le\sym(C)$ such that for any concrete category $\fK$ and any two Cayley objects $X,Y\in\fK$ over~$C$,
\qtnl{121117a}
\iso_\fK(X,Y)\ne\varnothing\quad\Leftrightarrow\quad
\iso_K(X,Y)\ne\varnothing.
\eqtn
\ethrm

The proof of Theorem~\ref{081117a} is given in Section~\ref{180418a}. The group $K$ constructed there is permutation isomorphic to the iterated wreath product 
$$
K=\AGL(1,p_1)\wr\cdots\wr\AGL(1,p_d),
$$
 where $p_1\ge\cdots\ge p_d$ are primes such that $n=p_1\cdots p_d$. One can replace the group~$K$ by a smaller group, e.g., the Hall $\pi$-subgroup of $K$, where $\pi=\{p_1,\ldots,p_d\}$. However,  it is doubtful that the order of such a group can be bounded from above by a polynomial in~$n$.\medskip
 
In principle, Theorem~\ref{081117a} could be used to test isomorphism of circulant combinatorial objects in time polynomial in their sizes. Indeed, the only thing  we need is to find a faithful and efficiently computable functor $F$ from the corresponding concrete category to the category of \{0,1\}-strings. If such a functor is given then to test isomorphism of the initial objects it suffices to check $F(K)$-isomorphism of the obtained strings and this can be done by the Babai-Luks algorithm \cite{BL83} in polynomial time, because the group~$K$ is solvable.\medskip 

For the concrete category of relational structures, considered in the paper, we use a different approach. We represent relational structures by special colored  hypergraphs in such a way that the required isomorphisms could be taken inside a solvable group constructed from the one mentioned in Theorem~\ref{081117a}. Finding these isomorphisms in polynomial time can be done with the help of the Miller's algorithm testing isomorphism of colored hypergraphs~\cite{Mil1983b} (see Section~\ref{250618a}). 

\thrml{211117a} 
The isomorphism of any two circulant objects can be tested in time polynomial in their sizes.
\ethrm 

It should be noted that Theorem~\ref{211117a}  cannot  be applied  directly to circulant hypergraphs. Indeed, to convert a hypergraph to a relational structure ina direct way, one should replace each hyperedge of cardinality $m$ with $m!$ tuples of length~$m$. But in this case the size of the resulting object may grow exponentially (of course, this is not the case if $m$ is a constant). Nevertheless, for circulant hypergraphs one can use the above mentioned Miller's algorithm and Theorem~\ref{081117a} to prove the following statement.

\thrml{211117a1} 
The isomorphism of any two circulant hypergraphs can be tested in time polynomial in their sizes.
\ethrm 

\medskip

All undefined notation and standard facts about permutation groups used in the paper can be found in the monographs~\cite{DM} and~\cite{W64}. 
In addition, we use the following notation.

\medskip

$\Omega$ denotes a finite set of cardinality $n$ and $\sym(\Omega)=\sym(n)$ is the symmetric group on $\Omega$.

The restriction of a group $K$ to a $K$-invariant set $\Delta\subseteq\Omega$ is denoted by $K^\Delta$.

The pointwise and setwise stabilizers of the set $\Delta$ in the group~$K$  are denoted by $K_\Delta$ and $K_{\{\Delta\}}$, respectively; we also set $K^\Delta=(K_{\{\Delta\}})^\Delta$.

For an imprimitivity system $\fD$ of a group $K\le\sym(\Omega)$, we denote by $K^\fD$  the permutation group induced by the action of $K$ on the blocks of~$\fD$.

\section{Proof of Theorem~\ref{081117a}}\label{180418a}

\subsection{}  
Let $C$ be a cyclic group of order $n=p_1\cdots p_d$, where $p_i,i=1,...,d$ are the prime factors of $n$. In what follows, we assume that for $i=1,\ldots, d-1$, 
$$ 
p_i\ge p_{i+1}\qaq n_i=p_1\cdots p_i.
$$ 
For any divisor $m$ of $n$, denote by $C_m$ a unique subgroup of~$C$ of order~$m$. Denote  by $\fD_m(C)$ the partition of $C$ into cosets of $C_m$ in~$C$. Clearly, $\fD_m(C)$ consists of $n/m$ classes, each of size~$m$.\medskip

Let $C_R$ be the (regular) subgroup of $\sym(C)$ induced by multiplications of~$C$. For any divisor $m$ of $n$, the partition $\fD_m(C)$ is an imprimitivity system  for~$C_R$. Its blocks coincide with the orbits of the group $(C_R)_m$. Note that although $(C_R)_m$ and $(C_m)_R$ are isomorphic as abstract groups they are distinct as permutation groups: the first one acts semiregularly on $C$ whereas the second one acts regularly on $C_m$.\medskip
 
In what follows we also set
$$
\Sup(C)=\{K\le\sym(C):\ C_R\le K\}.
$$
Thus, for any $K\in\Sup(C)$, the set 
$$
\cyc(K)=\{H\le K:\ H\ \,\text{is regular and cyclic}\}
$$
is not empty. A partial order $\preceq$ on $\Sup(C)$ is defined by the following rule: 
$$
L\preceq K\quad\Leftrightarrow\quad \forall H\in\cyc(K)\ \, \exists k\in K:\  \,
H^k\leq L.
$$ 
The following statement collects the results established in Theorems~1.8 and~4.9 of~\cite{Muz1999}.

\thrml{141117a} 
Any $\preceq$-minimal group $K\in\Sup(C)$ is solvable. Moreover, the partitions $\fD_{n_i}(C)$, $i=1,\ldots,d$, are imprimitivity systems for $K$.
\ethrm
 
\subsection{} 
Our aim is to construct a solvable subgroup in $\Sup(C)$ which is ``universal''  in the  sense that it contains every $\preceq$-minimal subgroup. To this end, let us define a group $\Wr(C)$ inductively by the number $d$ of primes in the decomposition of~$n=|C|$. Namely, set  $p=p_1$, $P=(C_R)_p$, $\fD_p=\fD_p(C)$, $\ov C=C/C_p$, and  
\qtnl{eq:Wr}
\Wr(C)=\{g\in\sym(C):\ \fD_p^g = \fD_p\, \land\, g^{\fD_p}\in \Wr(\ov C)\, \land\, [P,P^g]=1\}.
\eqtn
Note that in the case of $n=p$ we obtain $\Wr(C)=\bN_{\sym(C)}(C_R)\cong  \AGL_1(p)$, where $\AGL_1(p)$ is taken in its standard action on the set $\ov p=\{0,\ldots,p-1\}$, i.e., contains the $p$-cycle $(0,1,\ldots,p-1)$. 

\thrml{151118a} 
The following statements  hold:
\nmrt
\tm{1} $\Wr(C)$ is solvable,
\tm{2} $\fD_{n_i}(C)$ is an imprimitivity system for $\Wr(C)$, $i=1,\ldots,d$,
\tm{3} any solvable $S\in\Sup(C)$ preserving  the partitions $\fD_{n_i}(C)$, $i=1,\ldots,d$, is contained in $\Wr(C)$.
\enmrt 
\ethrm
\proof All parts of the statement are proved by induction on the number $d$ of prime divisors of $n$. Each part is trivial if $n=p$.  So, in what follows $n/p > 1$ and $d\geq 2$. \medskip

{\sf{ Part (1).}}   By definition, the partition $\fD_p$ is $\Wr(C)$-invariant and also
\qtnl{061218a}
\Wr(C)^{\fD_p} \leq \Wr(\ov C).
\eqtn
By the  induction hypothesis the group $\Wr(\ov C)$ is solvable. Therefore ${\Wr(C)}^{\fD_p}$ is solvable too. By  the Kaloujnine--Krasner embedding Theorem \cite[Theorem~2.6A]{DM} the group $\Wr(C)$ can be embedded to the wreath product $\Wr(C)^{C_p}\wr \Wr(\ov C)$. Thus it is sufficient to show that the group $\Wr(C)^{C_p}$ is solvable.\medskip 

We claim that 
$$
\Wr(C)^{C_p}\leq \bN_{\sym(C_p)}((C_p)_R)\cong \AGL_1(p).
$$
 Indeed, any $g\in \Wr(C)^{C_p}$ has a form $g = h^{C_p}$ for a suitable $h\in \Wr(C)_{\{C_p\}}$. By definition of $\Wr(C)$, $[P,P^h]=1$. The mapping $x\mapsto x^{C_p},x\in\Wr(C)_{\{C_p\}}$ is a group homomorphism. Therefore, 
$$
[P^{C_p},(P^h)^{C_p}]=1.
$$ 
Now taking into account that
$$
 P^{C_p}=(C_p)_R\qaq (P^h)^{C_p}= (P^{C_p})^{h^{C_p}} = (P^{C_p})^g,
$$
we arrive at the following implications:
$$
[P^{C_p},(P^h)^{C_p}]=1 \quad\Rightarrow\quad [(C_p)_R,((C_p)_R)^g]=1 \quad\Rightarrow\quad g\in \bN_{\sym(C_p)}((C_p)_R),
$$ 
as required.\medskip

{\sf{ Part (2).}} $\Wr(C)$-invariance of the partition $\fD_{n_1}(C)=\fD_p$ follows directly from formula~\eqref{eq:Wr}. Let $i>1$. By the induction hypothesis the partition $\fD_{n_i/p}(\ov C)$ is $\Wr(\ov C)$-invariant. In view of  inclusion~\eqref{061218a}, this implies that $\fD_{n_i/p}(\ov C)$ is also $\Wr(C)^{\fD_p}$-invariant.  Since $\fD_{n_i}(C)$ is the full preimage of $\fD_{n_i/p}(\ov C)$, it is $\Wr(C)$-invariant.\medskip

{\sf{Part (3)}.} Let $S\in\Sup(C)$ be a solvable group respecting  $\fD_{n_i}(C)$, $i=1,\ldots,d$.  Since $\fD_p$ is an $S$-invariant partititon, the induced group $S^{\fD_p}$ is solvable and fixes the partitions $\fD_{n_i/p}(\ov C)$. By the induction hypothesis, 
$$
S^{\fD_p}\leq \Wr(\ov C).
$$ 
Thus the first and the second conditions of~\eqref{eq:Wr} are satisfied. It remains to show that $[P,P^g]=1$  for each $g\in S$. \medskip

Indeed, the groups $P$ and $P^g$ have the same orbits, namely, the classes of the partition $\fD_p$. If  $K$ is  such a class, then $P^K$ and $ (P^g)^K$ are regular cyclic groups of degree $p$. In addition, they generate a solvable subgroup of $\sym(K)$. By Burnside's theorem  a Sylow $p$-subgroup of a solvable transitive permutation group of prime degree $p$ is normal. Therefore, $P^K = (P^g)^K$, and, consequently,  
$$
[P^K,(P^g)^K]=1.
$$
 Since this equality holds for every $K\in\fD_p$, we obtain $[P,P^g]=1$.\eprf\medskip

The statement below provides an exact description of the group $\Wr(C)$.

\lmml{201117c1}
The permutation group $\Wr(C)$ is permutation equivalent to the wreath product $\Wr(C_p)\wr \Wr(\ov C)$. Any two full cycles contained in $\Wr(C)$ are conjugate. In particular, all regular cyclic subgroups of $\Wr(C)$ are conjugate in $\Wr(C)$.
\elmm
\proof First we prove permutation equivalence. The partition $\fD_p$ is an  imprimitivity system for $\Wr(C)$. By the Kaloujnine--Krasner Theorem there exists a bijection $f:C\rightarrow C_p\times \fD_p$ such that 
$$
\Wr(C)^f\leq W_1\wr \Wr(\ov C),
$$ 
where $W_1 =\Wr(C)^{C_p}\leq\sym(C_p)$. The group~$W_1$ being solvable is contained in $\bN_{\sym(C_p)}((C_p)_R)=\Wr(C_p)$. Thus,
$\Wr(C)^f$ is a subgroup of $\Wr(C_p)\wr \Wr(\ov C)$, and
$$
\Wr(C)\leq (\Wr(C_p)\wr \Wr(\ov C))^{f^{-1}}.
$$
 The latter group is a solvable subgroup of $\sym(C)$, contains $C_R$, and satisfies the assumptions of Part~(3) of Theorem~\ref{151118a}. Therefore,  
 $$
 \Wr(C)\geq (\Wr(C_p)\wr \Wr(\ov C))^{f^{-1}},
 $$ 
 and, consequently, $\Wr(C) = (\Wr(C_p)\wr \Wr(\ov C))^{f^{-1}}$.\medskip
 
To prove the second statement we use induction on the number $d$ of prime factors of $n$. If $d=1$, then $\Wr(C)$ is permutation equiavlent to $\AGL_1(p)$ and the statement is true. The induction step follows  from \cite[Lemma~3.17]{1997a}.\eprf\medskip

Now Theorem~\ref{141117a} implies the following statement.

\crllrl{201117a}
Every $\preceq$-minimal subgroup of $\sym(C)$ is contained in $\Wr(C)$.
\ecrllr 

Applying  Lemma~\ref{201117c1}  inductively we conclude that $\Wr(C)$ is permutation equivalent to the wreath product
$$
\Wr(C_{p_1})\wr \Wr(C_{p_2})\wr\cdots\wr \Wr(C_{p_d})\cong \AGL(p_1)\wr\AGL(p_2)\wr \cdots  \wr \AGL(p_d).
$$

\subsection{Proof of Theorem~\ref{081117a}.}

Let $K=\Wr(C)$. For any concrete category $\fK$ the inclusion $\iso_\fK(X,Y)\supseteq \iso_K(X,Y)$ is obviously true for all $X,Y\in\fK$. Thus it suffices to verify the implication $\Rightarrow$ in formula~\eqref{121117a} only. To this end, let  $X,Y$ be Cayley objects over~$C$, i.e.,
$$
C_R\leq \aut_\fK(X)\qaq C_R\leq \aut_\fK(Y).
$$
Assume that $\iso_\fK(X,Y)\ne\varnothing$. Let $f\in\iso_\fK(X,Y)$. Then 
$$
H:=C_R^{f^{-1}}\leq \aut_\fK(X).
$$

Let $L$ be a $\preceq$-minimal subgroup of $\sym(C)$ contained in $\aut_\fK(X)$. Then there exists $u\in \aut_\fK(X)$ such that $H^u\le L$. However, $L\le K$ by Corollary~\ref{201117a}. Thus, 
$$
H^u\le L\le K.
$$ 
By Lemma~\ref{201117c1}, there exists $k\in K$ such that $(H^u)^k=C_R$. Thus, 
$$
(C_R)^{f^{-1}uk}=(H^u)^k=C_R.
$$
Consequently, the permutation  $f^{-1}uk$ lies in the normalizer of $C_R$ in $\sym(C)$. This normalizer is contained in $K$ by part (3) of Theorem~\ref{151118a} . Therefore, $f^{-1}u\in K$ and hence $u^{-1}f\in K$. Taking into account that $u^{-1}\in\aut_\fK(X)$, we obtain that
$$
X^{u^{-1}f}=(X^{u^{-1}})^f=X^f=Y.
$$
Thus,  $u^{-1}f\in\iso_K(X,Y)$ and so $\iso_K(X,Y)\ne\varnothing$.\eprf

\section{Proofs of Theorems~\ref{211117a} and~\ref{211117a1}}\label{250618a}
 
In \cite{Mil1983b},  the following problem was studied:\vspace{1mm}
 
 {\bf Hypergraph Isomorphism in a Coset.} {\it Given two edge-colored hypergraphs~$X$ and $Y$ with the same vertex set~$\Omega$, and a group $K\le\sym(\Omega)$, find the coset $\iso_K(X,Y)$}.\medskip
 
 Given $r\ge 2$, denote by $\Gamma_r$ the class of all finite groups whose composition factors are isomorphic to subgroups of $\sym(r)$. Then the algorithm proposed in~\cite{Mil1983b} solves the above problem in time polynomial in sizes of $X$ and $Y$, whenever $r$ is a constant and $K\in \Gamma_r$. In fact, the only essential property of $\Gamma_r$ used for estimation the complexity of the algorithm is that the order of any primitive group of degree $n$ in $\Gamma_r$ is at most $n^{O(1)}$ for a fixed~$r$. In particular, the algorithm is still polynomial-time if the group $K$ is solvable, because the order of any solvable primitive group of degree $n$ is at most $n^4$ \cite{P82}.\medskip
 
 {\bf Proof of Theorem~\ref{211117a1}.}  Let $X$ and $Y$ be Cayley hypergraphs over a cyclic group~$C$. By Theorem~\ref{081117a}, the equivalence~\eqref{081117a}  holds for 
the category of hypergraphs. Since the group~$K$ is solvable, the Miller's algorithm tests whether the set $\iso_K(X,Y)$ is empty in time polynomial of the size of~$X$.\eprf\medskip

{\bf Proof of Theorem~\ref{211117a}.} We need an auxiliary construction associating a Cayley object $X$ over a cyclic group~$C$ with an  edge-colored hypergraph $\cH=\cH(X)$. The vertex set of $\cH$ is  the disjoint union 
 $$
\Omega=\Omega_1\cup\cdots\cup\Omega_m,
$$
where $\Omega_i$ is a copy of~$C$ for all~$i=1,\ldots,m$, and $m$  is the maximal arity of the relation entering the set 
$$
\cR=\cR(X)=\{R_1,\ldots,R_k\}
$$ 
of the relations defining~$X$ (the indices of the  $R_i$ are defined accordingly the linear ordering of~$\cR$). \medskip
 
 The hyperedge set of $\cH$ consists of  the three parts $E_1$, $E_2$, and $E_3$. Namely,
 $$
 E_1=\{\Omega_1,\ldots,\Omega_m\},
 $$ 
and the color of $\Omega_i$ is set to be~$i=1,\ldots,m$.  The second part $E_2$ is the union of the sets 
 $$
\wh R_i= \{\{\alpha_1,\beta_2,\ldots\}:\ (\alpha,\beta,\ldots)\in R_i\}, \quad i=1,\ldots,k,
 $$
where $\alpha_1,\beta_2,\ldots$ are the copies of the vertices $\alpha,\beta,\ldots$ of~$X$  belonging to $\Omega_1,\Omega_2,\ldots$, respectively; the color of any hyperedge from~$\wh R_i$ is defined to be $m+i$. Finally, the third part consists of $n$ hyperedges each of color $m+k+1$, 
$$
E_3=\{\{\alpha_1,\ldots,\alpha_m\}:\ \alpha\in G\}. 
$$
  
\lmml{250618u}
Given a permutation $f\in\sym(C)$, denote by $\wh f$ the permutation on~$\Omega$ defined by the formula
$$
(\alpha_i)^{\wh f}=(\alpha^f)_i,\qquad \alpha\in C,\ 1\le i\le m.
$$
Then for any Cayley object $Y$ over $C$,
$$
\iso(\cH(X),\cH(Y))=\{\wh f:\ f\in\iso(X,Y)\}.
$$
\elmm  
 \proof From the definition of the permutation $\wh f$, it follows that
 $$
(E_1)^{\wh f}=\Omega_i\qaq (E_3)^{\wh f}=E_3.
 $$
Assume that $f\in \iso(X,Y)$. Then  $\cR(Y)=\{(R_i)^f:\ 1\le i\le k\}$.  Again the definition of $\wh f$ implies that
$$
(\wh R_i)^{\wh f}=\wh{(R_i)^f}, \quad 1\le i\le k.
$$
Thus, $\wh f\in \iso(\cH(X),\cH(Y))$. The converse inclusion is verified in a similar way.\eprf\medskip

Now let the group $C$ be cyclic, $\fK$ the category of circulant objects, and $X,Y\in\fK$. Then to test isomorphism between $X$ and $Y$, it suffices to check that 
$$
\iso_K(X,Y)\ne\varnothing,
$$ 
see Theorem~\ref {081117a}.  However, from Lemma~\ref{250618u}, it follows   that
$$
\iso_K(X,Y)\ne\varnothing\quad\Leftrightarrow\quad \iso_{\wh K}(\cH(X),\cH(Y))\ne\varnothing,
$$
 where $\wh K=\{\wh k:\ k\in K\}$ is a permutation group on~$\Omega$ isomorphic to~$K$. Since the latter, and hence $\wh K$,  is solvable, the set $\iso_{\wh K}(\cH(X),\cH(Y))$ can be found by Miller's algorithm in time $\poly(nm)$. This completes the proof, because $m\le n$.\eprf

\end{document}